\documentclass{article}
\usepackage{amsmath,amsfonts,amssymb,theorem,verbatim}

\renewcommand{\abstractname}{Abstract.}
\renewcommand\abstract{\hfil\break\topsep=0pt%
\partopsep=0pt\parsep=0pt\itemsep=0pt\relax
\trivlist\item[\hskip\labelsep
{\bfseries\abstractname}]\if!\abstractname!\hskip-\labelsep\fi}

\newcommand{\email}[1]{{(e-mail: #1)}}

\def\keywordname{{\bfseries Key words:}}
\def\keywords#1{\par\addvspace\baselineskip\noindent\keywordname\enspace
\ignorespaces#1}

\newtheorem{Thm}{Theorem}
\newtheorem{Cor}{Corollary} \newtheorem{Lemma}{Lemma}
\newtheorem{Prop}{Proposition} 
{\theorembodyfont{\normalfont}
\newtheorem{Ex}{Example}

}

\def\BEN{\begin{enumerate}}  \def\BI{\begin{itemize}}
\def\EEN{\end{enumerate}}   \def\EI{\end{itemize}} \def\im{\item}
\def\ssec{\subsection}    
  \def\no{\noindent}

\def\nn{\nonumber}
\def\beq{\begin{eqnarray}} \def\eeq{\end{eqnarray}}

\def\eqn#1{\begin{equation}#1\end{equation}}

\def\aln#1{\begin{align}#1\end{align}}
\def\al*#1{\begin{align*}#1\end{align*}}

\def\ga*#1{\begin{gather*}#1\end{gather*}}

\def\alat*#1#2{\begin{alignat*}{#1}#2\end{alignat*}}
\def\bea{\begin{eqnarray*}} 
\def\eea{\end{eqnarray*}}
\def\ml*#1{\begin{multline*}#1\end{multline*}}

 \def\mbf{\mathbf} \def\mrm{\mathrm} 
  
\newcommand{\Bf}[1]{{\mbox{\scriptsize\boldmath$#1$}}}
\newcommand{\bff}[1]{{\mbox{\boldmath$#1$}}}

\newcommand{\Prob}{{\rm I\hspace{-0.7mm}P}}
\newcommand{\Exp}{{\rm I\hspace{-0.7mm}E}}

\def\P{{\mathbb P}} 
  \def\R{{\mathbb R}}

\def\mc{\mathcal}

\def\le{\left} \def\ri{\right} \def\i{\infty}

\def\te#1{\mathrm{e}^{#1}}  \def\td{\text{\rm d}}
 
\def\I{\int}     

  \def\H{\hat}
\def\WH{\widehat} \def\WT{\widetilde}

\def\a{\alpha} \def\b{\beta}
\def\g{\gamma}  \def\d{\delta}    

\def\e{\epsilon} \def\k{\kappa} \def\l{\lambda} \def\m{\mu} \def\n{\nu}
  \def\nn{\nonumber}  \def\r{\rho} \def\s{\sigma}
   \def\f{\varphi}  \def\ps{\psi} 
  \def\q{\qquad} \def\D{\Delta} 
\def\F{\Phi} \def\G{\Gamma}

\newcommand{\proof}{\no{\it Proof\ }}

\newcommand{\exit}{{\mbox{\, \vspace{3mm}}} \hfill\mbox{$\square$}}
\begin{document} 
\author{Martijn 
Pistorius\footnote{Department of Mathematics, Strand, London WC2R 2LS, UK,
{\tt\email{pistoriu@mth.kcl.ac.uk}}}}
\title{On maxima and ladder processes
for a dense class of L\'{e}vy processes}

\date{\sl King's College London}
\maketitle

\begin{abstract}
Consider the problem to explicitly calculate the law of 
the first passage time $T(a)$ of a general L\'{e}vy process $Z$
above a positive level $a$. In this paper 
it is shown that the law of $T(a)$
can be approximated arbitrarily closely by the laws of $T^n(a)$,
the corresponding first passages time for $X^n$, where 
$(X^n)_n$ is a sequence of  L\'{e}vy processes
whose positive jumps follow a {\it phase-type} distribution.
Subsequently, explicit expressions are derived for the laws of $T^n(a)$ 
and the upward ladder process of $X^n$.
The derivation is based on an embedding of $X^n$ into
a class of Markov additive processes and on the solution 
of the fundamental (matrix) Wiener-Hopf factorisation 
for this class.
This Wiener-Hopf factorisation 
can be computed explicitly by solving iteratively
a certain fixed point equation. It is shown that, 
typically, this iteration converges geometrically fast.
\end{abstract}

\smallskip

\no

\section{Introduction}
A L\'{e}vy process is a stochastic process with c\`{a}dl\`{a}g paths and 
homogeneous independent increments. L\'{e}vy processes have turned up 
in a number of areas of applied probability, 
for instance as models for the workload of queues, 
the content of a dam and the reserve level of an insurance company, 
see e.g. \cite{aruin,aque,prabhu} and references therein.
More recently, L\'{e}vy processes have also been proposed as models for 
the evolution of the logarithms of interest rates or prices of assets 
and as model for credit derivatives, see e.g. \cite{CarSch,Eberlein,MCC} 
to name just three articles.
In many of the mentioned areas one is interested in an 
(explicit) characterisation of the distributions of the running maximum, 
the first passage time over a level and the {\em ladder height} 
and {\it ladder time} processes.

The study of pathwise extrema of a stochastic
process is called fluctuation theory. 
For L\'{e}vy processes with only jumps of one sign, 
the fluctuation theory is simplified and results are 
more explicit than in the general case; 
See e.g. the comprehensive review of Bingham 
\cite{bingham}, Chapter VII in Bertoin \cite{bert96} 
or the reviews from martingale \cite{KPa} 
or potential-theoretic perspective \cite{ptnt}. 
In the case of a general L\'{e}vy process with positive {\it and} 
negative jumps few explicit results are available.  
However, restricting one-self to the class of 
L\'{e}vy process with positive jumps of {\it phase-type} 
(see Section \ref{phtd} below for definition),
Asmussen et al. \cite{AAP} explicitly characterised 
the joint law of the first 
exit time from the negative half-line and the corresponding 
overshoot. Using this result, the law of the first passage time 
of a positive level for a general L\'{e}vy process can be 
approximated arbitrarily closely, as shown in Section \ref{ssec:phtlevy}.
Indeed, any L\'{e}vy process can be approximated 
arbitrarily closely by a L\'{e}vy process with phase-type jumps  
and this convergence is shown to carry over to the first exit times mentioned. 

In \cite{AAP} martingale techniques were invoked to characterise 
the joint law of the aforementioned stopping time and overshoot as 
solution of a certain linear system.
In this paper, we study the same class of L\'{e}vy processes $X$
but follow a different approach from that of \cite{AAP}.
First, the L\'{e}vy process is embedded into a class of Markov additive processes 
(or a phase-process perturbed by a spectrally negative L\'{e}vy process) 
and then the matrix Wiener-Hopf factorisation is characterised. 
This factorisation extends existing results in the literature: 
We mention in particular Asmussen \cite{a} and Rogers \cite{rogers}, 
who covered the factorisation of random walk
and phase processes perturbed by Brownian motion, respectively.

The analysis of this Wiener-Hopf factorisation 
leads to an alternative characterisation of the law 
of the maximum of $X$ and to an explicit 
description of the law of the upcrossing ladder process of $X$.
For explicit calculation of this factorisation, an algorithm is derived and 
it is shown that, typically, this algorithm 
converges geometrically fast.

The rest of the paper is organised as follows.
In Section \ref{sec:prel} the notation
is set and some theory regarding the Wiener-Hopf 
factorisation and ladder processes of  
L\'{e}vy processes is revisited. 
Section \ref{sec:MWE} is then devoted to the matrix Wiener-Hopf
factorisation of the process $X$ and its embedding. 
In Sections \ref{sec:fluc} explicit identities
for the law of the first passage time over a positive level and
the law of the ladder process are derived. 
Finally, in Section \ref{sec:iter}, 
an algorithm for the computation of the aforementioned laws 
and its convergence are studied.

\section{Preliminaries}\label{sec:prel}
\ssec{Phase-type distributions}\label{phtd}
A distribution $F$ on $(0,\infty)$ is {\it phase-type} if it is
the distribution of the absorbtion time $\zeta$ in a finite state
continuous time Markov process $J=\{J_t\}_{t\ge 0}$ with one state
$\Delta$ absorbing and the remaining ones $1,\ldots,m$ transient.
That is, $F(t)=$ $\P (\zeta\leq t)$ where $\zeta=$ $\inf\{s>0:\,
J_s=\Delta\}$. The parameters are $m$, the restriction $\bff{T}$
of the full intensity matrix to the $m$ transient states and the
initial probability (row) vector $\bff{\alpha}=$
$(\alpha_1\,\ldots\,\alpha_m)$ where $\alpha_i=\P(J_0=i)$. For any
$i=1,\ldots,m$, let $t_i$ be the intensity of a transition
$i\to\Delta$ and write $\bff{t}=$ $(t_1\,\ldots\,t_m)'$ for the
(column) vector of such intensities. Note that $\bff{t}= -\bff{T}
\bff 1,$ where $\bff 1$ denotes  a column vector of ones. It
follows that the cumulative distribution $F$ is given by:
\eqn{\label{phtDF}
 1- F(x)=\bff{\alpha}\te{\Bf{T}x} \bff 1,
} the density is $f(x)=$ $\bff{\alpha}\te{\Bf{T}x} \bff t$ and
the Laplace transform is given by
$
\hat{F}[s]=\I_0^\i\te{-sx}F(\td x) =
\bff{\alpha}(s\bff{I}-\bff{T})^{-1}\bff{t}.
$
Note that
$\hat{F}[s]$ can be extended to the complex plane except at a
finite number of poles (the eigenvalues of $\bff{T}$).
{ 
Throughout we will assume that the 
representation of the form (\ref{phtDF}) for the distribution function $F$ 
is {\em minimal}, that is, there exists no number $k<m$, $k$-vector $\bff b$ 
and $k\times k$-matrix $\bff G$ such that $1-F(x) = \bff b\te{\Bf G x}\bff 1$.
}
\subsection{Phase--type L\'{e}vy processes}\label{ssec:phtlevy}
Consider  a L\'{e}vy process $X$ of the form
\begin{equation}\label{Levygenal}
X_t = X^{(+)}_t + X^{(-)}_t,
\end{equation}
where $X^{(-)}=\{X^{(-)}_t, t\ge0\}$ is a L\'{e}vy process
without positive jumps and
$X^{(+)}=\{X^{(+)}_t, t\ge0\}$
is an independent compound Poisson process with intensity
$\l^{(+)}$
and jumps of phase-type with parameters
$(m^{(+)}, \bff T^{(+)}, \bff\a^{(+)})$. We exclude the case 
that $X^{(-)}$ is a negative deterministic drift.
We define by $\k(s) = \k_X(s) = \log\Exp[\te{sX}]$
the L\'{e}vy exponent of $X$ which is well defined at least 
for $s$ with $\Re(s)=0$ and which is in this case given by
$$
\k(s) = \k_{X^{(-)}}(s) + \l^{(+)}(\hat F^{(+)}[-s]-1),
$$
where $\k_{X^{(-)}}$ denotes the L\'{e}vy exponent of $X^{(-)}$ 
and $\hat F^{(+)}[s] = \bff\a^{(+)}(s\bff I - \bff T^{(+)})^{-1}\bff t$
is the Laplace transform of $F^{(+)}$. Note that $\k$ can be analytically 
extended to the positive half-plane except a finite number of poles 
(the eigenvalues of $-\bff T^{(+)}$) and we shall denote 
this extension also by $\k$.

Any L\'{e}vy process $L$ may be approximated arbitrarily closely
in law by a sequence $(X(n))_n$ of L\'{e}vy processes of type (\ref{Levygenal}).
Indeed, it is well known that $X(n)$ weakly converges to $L$ (as a process, in the
Skorokhod topology) if and only if $(X(n)_1)_n$ converges 
to $L_1$ in distribution (see e.g. Jacod and Shiryaev \cite{JS}, 
Cor. VII.3.6). Also, the set of phase-type distributions is dense 
in the set of probability distributions on $(0,\i)$
(in the sense of weak convergence of 
probability distributions). 
Therefore the aforementioned approximation can be obtained 
by first restricting the L\'{e}vy
measure $\nu$ of $L_1$ to $\R\backslash [-\e,\e]$ and then using 
that the probability distributions $c^+_\e\mbf 1_{\{x>\e\}}\nu(\td
x)$ and $c^-_\e \mbf 1_{\{x<-\e\}}\nu(\td x)$ (with $c^+_\e=1/\nu(\e,\i)$
and $c^-_\e=1/\nu(-\i,-\e)$ and $\mbf 1_A$ the indicator of the set $A$) 
can be approximated arbitrarily closely
by phase-type distributions. The relevant
methodology for fitting a phase--type distribution to a given
 set of data is developed in \cite{ANO} for
traditional maximum likelihood and in \cite{BGL} in a Bayesian
setting.

Write $T^+(a) = T^+(a)(X)$ for the first passage time of $X$ over $a$
\begin{equation}\label{eq:fipa}
T^+(a) = \inf\{t\ge 0: X_t > a\}
\end{equation}
and denote by $O^+(a)(X)=X_{T^+(a)} - a$ the corresponding overshoot of $X$. 
The next result shows that the weak convergence of the processes 
$X(n)$ carries over to the stopping times $T^+(a)(X(n))$ and 
the overshoots $O^+(a)(X(n))$
\begin{Prop}\label{prop:convstopp}
Let $(X(n))_n$ be a sequence of L\'{e}vy processes such that
$X(n) \to X$ weakly in the Skorokhod topology. Then, as $n\to\i$,
$$
(T^+(a)(X(n)), O^+(a)(X(n)) \longrightarrow (T^+(a)(X), O^+(a)(X))
$$
where the convergence is in distribution.
\end{Prop}
Before giving the proof we first review Wiener-Hopf factorisation 
of L\'{e}vy processes.
\subsection{Wiener-Hopf factorisations} In this subsection, we review
some of the fluctuation theory of L\'{e}vy processes. For more
background we refer the reader to Bingham \cite{bingham} or
Bertoin \cite{bert96}. Denote by $\mc I^{(+)}=\{i: \Re(\r_i) >
0\}$ the set of roots $\r_i$ with positive real part of the
Cram\`er-Lundberg equation \eqn{\label{CL} \k(\r) = \k_X(\r) = a.
} Let $x\wedge 0 = \min\{x,0\}$ and $x\vee 0 = \max\{x,0\}$ 
and write $M_t=\sup_{s\leq
t}(X_t\vee 0)$ and $I_t=\inf_{s\leq
t}(X_t\wedge 0)$ for the supremum and infimum of $X$ up to time $t$, 
respectively. Similarly, write $M^{(-)}_t$ and $I^{(-)}_t$ for 
the the corresponding quantities of $X^{(-)}$.
Denote by $e(a)$  an independent exponential random variable
with mean $a^{-1}$. Set for $s$ with $\Re(s)\geq 0$
\eqn{\label{eq.wienerhopf2} \f_a^-(s) = \Exp[\exp(sI_{e(a)})],\q
\f_a^+(-s) = \Exp[\exp(-sM_{e(a)})]. } The functions
$s\mapsto\f^\mp_a(s)$ are analytic for $s$ with $\pm\Re(s)>0$,
respectively. By bounded convergence it follows that
$\f_a^-(\i)=\Prob(I_{e(a)}=0)$ and
$\f_a^+(-\i)=\Prob(M_{e(a)}=0)$. For $a>0$, the functions
$s\mapsto\f^\mp_a(s)$ satisfy the Wiener-Hopf factorisation
\eqn{\label{eq.wienerhopf0} a/(a-\k(s))=\f^+_a(s)\f^-_a(s)
\q\text{for all $s$ with $\Re(s)=0$} } For a proof we refer to e.g. 
Bingham \cite[Thm. 1a]{bingham}. 
Since $|\f^+_a(s)| = |\f^-_a(s)| \leq 1$ for $s$
with $\Re(s)=0$, there are no roots of (\ref{CL}) with zero real
part when $a>0$. 
If $X$ is a L\'{e}vy process of 
the form (\ref{Levygenal}), the phase-type representation is minimal 
and $a>0$, then 
$\f^+_a$ is explicitly given by (as shown in \cite{AAP}) 
\eqn{
\label{eq.wienerhopf}
\f_a^{+}(s) = \frac{\det(-s\bff I-\bff T^{(+)})}{\det(-\bff T^{(+)})}\cdot
\frac{\prod_{i\in\mc I^{{(+)}}}(-\r_i)}
{\prod_{i\in\mc I^{(+)}}(s-\r_i)}, 
}
where the first factor is to be taken equal to 1 if $X$ has no negative 
jumps.
The following Wiener-Hopf identity (e.g. \cite[Thm. 1e]{bingham}) 
links the joint distribution 
of the first time of crossing the level $a$ 
and the corresponding overshoot $(T^+(a), O^+(a))$ 
to the Wiener-Hopf factor $\f_a^+$:
\begin{equation}\label{bing-thm1e}
\Exp\le[\te{-qT^+(e(\l)) - \m O^+_{e(\l)}}\ri]
= \frac{\l}{\l-\m}\le[1-\frac{\f_q^+(-\l)}{\f_q^+(-\m)}\ri].
\end{equation}
We have made all preparations for the proof of Proposition \ref{prop:convstopp}.
\medskip

\proof{\it of Proposition \ref{prop:convstopp}} 
Writing $M^n_t = \sup_{0\leq s \leq t}(X(n)_s\vee 0)$ 
for the running supremum of $X(n)$, the triangle 
inequality implies that
$|M^n_t- M_t|$ is smaller than $\sup_{s\leq t}|X(n)_s\vee 0 -
X_s\vee 0|$ which is smaller than $\sup_{s\leq t}|X(n)_s - X_s|$. 
Since, by assumption, $X(n)$ converges to $X$ in the Skorokhod topology 
and since a L\'{e}vy process is continuous at each fixed times a.s., 
it can be directly verified (as in Prop. VI.2.4 in \cite{JS}) that
$M^n_t$ converges in distribution to $M_t$ for fixed $t\ge0$. 
This implies that the
moment-generating function (mgf) of $M^n_{e(q)}$ converges to the
mgf of $M_{e(q)}$. By the Wiener-Hopf identity \ref{bing-thm1e}
and the extended continuity theorem it follows then that the joint
Laplace transforms of $(T^+(a)(X(n)), O^+(a)(X(n)))$
converges to the joint Laplace transform of $(T^+(a)(X), O^+(a)(X))$, 
which shows the stated convergence in distribution.\exit
\medskip

Closely related to the supremum process $M$ are the ascending {\em
ladder time} and the {\em ladder height} processes, which are 
L\'{e}vy processes. One of the objectives of this note is to
explicitly find the distribution of these two processes. To be
able to describe the ladder process, we first need look at the
{\em local time} of $M-X$ at zero. We will distinguish between two
different cases (see \cite[Ch. IV]{bert96} for details):

If $X^{(-)}$ is {\it not} the negative of a subordinator, $0$ is
regular for $M-X$ and the local time process can be taken to be
continuous. The canonical choice in this case is to take $L =
M^c$, the continuous part of the supremum process $M$.

If $X^{(-)}$ is the negative of a subordinator, $0$ is irregular
for $M-X$ (i.e. $\Prob_0(\s=0)=0$, where $\s$ is the first time
$M-X$ reaches 0) and the zero set of $M-X$ forms a discrete set.
The supremum is a jump process, where the jump sizes form an
i.i.d. sequence and the jump times are precisely this zero set. In
order to let it fit in the same framework as the previous case, an
extra randomisation is needed, which can be thought of as the
analog of the random time change to transform a random walk into a
compound Poisson process. Let $e_1(c), e_2(c), e_3(c), \ldots$ be
a sequence of i.i.d. exponential random variables with a certain
parameter $c>0$ and denote by $n(t)$ the number of zeros of $M-X$
up to but excluding time $t$, i.e. $n(t) = \max\{i : R_i < t\}$
where $R_i$ are the subsequent zeros of the process $M-X$. Then
the local time $L_t$ at time $t$ is given by
\eqn{\label{defLs}
L_t = \sum_{n=1}^{n(t)}e_n(c).
}
The ladder time is the right-inverse $L^{-1}$
of the local time $L$
$$
L^{-1}_t = \inf\{s\ge0: L_s > t\}
$$
and is a stopping time with respect to the standard filtration
generated by $X$. The ladder height process $H$ is taken to be $H
= M_{L^{-1}}$, the supremum at this stopping time. Write $\k^+$
for the joint characteristic exponent of the ladder process
$(L^{-1},H)$
$$
\exp\{ - \k^+(\a,\b) \}= \Exp\le[\exp\{ - \a L^{-1}(1) - \b H(1)\}\ri]
,\q \a,\b >0.
$$
Denote by $G_t=\sup\{u\leq t: X_u = M_u\}$ the last time before $t$ that $X$
was at its supremum and let as before $e(q)$ be an independent
exponential time with mean $q^{-1}$. Then an extension of 
the basic Wiener-Hopf factorisation (see e.g. \cite[Ch
VI]{bert96} or \cite{bingham}) tells us that 
$(G_{e(q)}, M_{e(q)})$ and
$(e(q)-G_{e(q)}, M_{e(q)} - X_{e(q)})$ are independent and the joint distribution 
of $(G_{e(q)}, M_{e(q)})$ is expressed in terms of $\k^+$ by 
\begin{equation}\label{WH3}
\Exp[\te{-a G_{e(q)} - bM_{e(q)}}] =
\frac{\k^+(q,0)}{\k^+(q+a,b)}\q a,b>0.
\end{equation}
Note that this expression and the stated independence imply that
\begin{equation}\label{WH4}
\Exp[\te{-a (e(q)-G_{e(q)}) - b(M_{e(q)}-X_{e(q)})}] =
\frac{q}{\k^+(q,0)}\frac{\k^+(q+a,b)}{q+a - \k(b)}.
\end{equation}
Similarly, we can define $\k^-$ to be the joint L\'{e}vy exponent
of $(\WH L^{-1},\WH H)$, the ladder process of the dual $\WH X=-X$
of $X$ and then we have the following relation between the different
characteristic exponents $\k, \k^+, \f^+_a$ for $s$ with $\Re(s)=0$ and $a>0$
(e.g. Bertoin \cite[Ch. VI]{bert96}):
\aln{\label{WH2}
\f^+_a(s) &= \frac{\k^+(a,0)}{\k^+(a,-s)},\\
\k(s) &= -C\cdot\k^+(0,-s)\k^-(0,s),\label{WHkappa}
}
where $C$ is some constant. 
The second identity in (\ref{WHkappa}) 
is often referred to as the Wiener-Hopf factorisation of the
L\'{e}vy exponent $\k$.
\subsection{Embedding}\label{sec:emb}
To study properties of the L\'{e}vy process $X$, we follow ideas 
of \cite{AAP} and embed $X$ into a Markov process $(A, Y)$.
Informally, ones get $A$ from $X$ by  'levelling out' the positive jumps into 
piecewise linear parts of gradient $+1$; the process $Y$ is set equal to zero
if $X$ moves like $X^{(-)}$ and equal to the underlying phase process of 
an upward jump otherwise.

More precisely and slightly more generally, 
let $Y$ be an irreducible Markov process with finite state space 
$E$, with $E = \{0,1,\ldots, m\}$. 
Denote the generator of $Y$ restricted to $E$ by $\bff Q$. 
Letting $X^{(-)}$ be the spectrally negative L\'{e}vy process
of (\ref{Levygenal}) and $v$ and $\s$  
functions that map $E$ to $\R$,  
the additive functional $A=\{A_t\}_{t\ge0}$ is defined as
\eqn{\label{eq:AXmin}
A_t = A_0 + \I_0^t\s(Y_s)\td X^{(-)}_s + \I_0^tv(Y_s)\td s.
}
The pure fluid model corresponds to setting 
$\s$ equal to zero whereas this model reduces to a 
spectrally one-sided L\'{e}vy process if $v\equiv0$, $\s$ is 
constant and $Y$ is recurrent.

Next, we define $\WT Y^+$ to be the {\it upcrossing ladder process} 
of $Y$, that is, 
\eqn{
\label{eq:jtilde}
\WT Y^+_t = Y(\g_t) \text{ where } \g_t = \inf\{s\ge0: A_s > t\}.
}
It is easily verified that this time-changed process $\WT Y$ 
is again a Markov process; we denote its generator by $\bff Q^+$.
The next section is devoted to a characterisation of its form.

Consider now the special case of above additive functional 
(\ref{eq:AXmin}) where $m=m^{(+)}$, 
and restriction of the intensity matrix of $Y$ to $\{0,1,\ldots, m\}$ 
is given by (in block notation) $\bff Q_0$, where for $a\ge0$
\eqn{\label{Qa}
\bff Q_a = \begin{pmatrix}
-\l^{(+)}-a& \l^{(+)}\bff\a^{(+)}\\
\bff t^{(+)} & \bff T^{(+)}
\end{pmatrix}.
}
Moreover, for $i\in E$ set $v(i) = 1-\s(i)$ and $\s(i) = \d_{0i}$, the Kronecker delta,  and 
let $T_0(t)=\I_0^tI(Y_s=0)\td s$ denote local time of $Y$ at $0$, that is, 
the amount of time before time $t$ that $Y$ has spent in state $0$. Then we get 
back the original process $X$ by time-changing $A$ with the inverse local time
\begin{equation}
T_0^{-1}(u) = \inf\{t\ge0: T_0(t) > u\}. 
\end{equation}
Indeed, the independence of the increments of $X$ implies that $X$ is in law equal to 
$A \circ T_0^{-1}$. 
Exponential killing of the original L\'{e}vy process $X$ at rate $a$ 
can be incorporated by replacing $\bff Q_0$ by $\bff Q_a$, since all states of $Y$ 
other than $0$ originate from positive jumps of $X$ so that the local time $T_0$ of
$A$ at zero corresponds to the time scale of $X$.

\section{Matrix Wiener-Hopf factorisations}\label{sec:MWE}
Denoting by $v$ and $\s$ functions that map $E$ to $\R$ and $[0,\i)$ 
respectively, we consider now the additive functional $A=\{A_t, t\ge0\}$ 
given by (\ref{eq:AXmin}).
In the sequel, we restrict ourselves to a characterisation
of the generator $\bff Q^+$ of the upcrossing ladder process 
$\WT Y^+$ of $A$ given in (\ref{eq:AXmin}). In the analysis, 
we shall distinguish between the cases that $-X^{(-)}$ is a subordinator, 
(that is, $X^{(-)}$ has non-increasing paths) or that $X^{(-)}$ 
has non-monotone paths. 
We partition the state space $E$ of $Y$ 
into a part $E^-$, where $v(j)\leq 0$, and a part $E^+$,  where 
$v(j)>0$. It follows that $A$ decreases
as long as $Y$ is in one of the  states in $E^-$ and that $A$ 
can decreases as well as increase when $Y$ is $E^+$ and
the state space of the upcrossing ladder process 
$\WT Y^+$ defined in (\ref{eq:jtilde}) is given by $E^+$. 
The matrix of ``up-crossing phase probabilities'' is denoted 
by $\bff\eta$, that is,
\eqn{\label{etaij} \bff\eta(i,j) = \Prob[\WT Y^+_0=j|Y_0 = i],\q
i\in E^-,\q j\in E^+. }
Finally, let $\bff K(\s, -\bff G)$ be the matrix whose rows are given by
$(\bff 1_i\kappa_{X^{(-)}}(-\s(i)\bff G))$ with $\kappa_{X^{(-)}}$
given by (by the L\'{e}vy-Khintchine formula)
\eqn{
\label{eq:kappamin}
\k_{X^{(-)}}(-\bff G)=
\frac{s^2}{2}\bff G^2 - c \bff G + \I^0_{-\i}\le(\te{-\Bf G x} -
\bff I -x\bff G\mbf 1_{\{|x|<1\}}
\ri)\nu(\td x),
}
with $\n$ the L\'{e}vy measure and $s$ the Gaussian coefficient of $X^{(-)}$.
\begin{Thm}\label{uni-s}
(i) If $-X^{(-)}$ is a subordinator,
the generator matrix\ $\bff Q^+$ of the process
$\WT Y^+$ satisfies the matrix equation given, in block notation, by
\begin{equation}
\label{systemnoiseWHs}
\begin{pmatrix}\bff\eta\\\bff I\end{pmatrix}
\bff K(\s, -\bff G)
+ \bff Q \begin{pmatrix}\bff\eta\\\bff I\end{pmatrix}
= \bff V\begin{pmatrix}\bff\eta\\\bff I\end{pmatrix}\bff G,
\end{equation}
(ii) {If $\bff Q$ is transient or if $\bff Q$ is recurrent and $\sup_t A_t=\i$ a.s., 
the solution in $\bff G\in\mathcal Q(|E^+|)$ of (\ref{systemnoiseWHs}) is unique.}
\end{Thm}

Now we turn to the case that $A$ is given by 
(\ref{eq:AXmin}), where $-X^{(-)}$ is not a subordinator. 
We restrict ourselves to the case that the function $v$ is positive, 
$v(i)>0$ , for each $i\in E$ for which $\s(i)=0$, 
to ensure that $A$ can increase in each state in $E$.

\begin{Thm}\label{uni-ns}
(i) If $-X^{(-)}$ is not a subordinator,
the generator matrix\ $\bff Q^+$ of the process
$\WT Y^+$ defined in (\ref{eq:jtilde}) solves the equation
\begin{equation}
\label{systemnoiseWH}
\bff K(\s,-\bff G) + \bff Q = \bff V\bff G.
\end{equation}

(ii) If $\bff Q$ is transient or if $\bff Q$ is recurrent and $sup_t A_t = +\i$
a.s., $\bff Q^+$ is the unique $\bff G\in\mathcal Q(|E|)$ that
solves  (\ref{systemnoiseWH}).
\end{Thm}

\proof{\it of Theorem \ref{uni-ns}}
(i)
Define the function $f: E\times \R\to \R$ by
$
f(j,x) = \Exp_{j,x}[h(\WT Y^+_k)\mbf 1_{\{\g_k < \i\}}]
$
where $h$ is any function on $E$
and where $\Exp_{j,x}$ denotes the expectation under the measure
$\Prob$ conditioned on $\{Y_0=j, A_0=x\}$.
Since $\WT Y^+$ is a Markov process with generator $\bff Q^+$, 
the function $f$ can be expressed as
\eqn{\label{eq:mkg}
f(i, x) = \bff 1_i'\exp(\bff Q^+(k-x))\bff h,\q i\in E,\q x\leq k.
}
However, on the other hand, it is straightforward to check that
$f(Y_{t\wedge\g_k}, A_{t\wedge\g_k})$
is a martingale and we find by It\^o's
lemma that $\bff f = (f(i,u),i\in E)$ satisfies
for $i\in E$ and $u < k$,
\begin{equation}
\label{ITOPHASE}
\G(\s(i)f(i,u)) + v(i)f'(i,u) + \sum_j q_{ij}(f(j,u)-f(i,u)) = 0, 
\end{equation}
where $\G$ denotes the infinitesimal generator of the process $X^{(-)}$:
$$
\G f(x) = \frac{\s^2}{2}f''(x) + cf'(x) + \I_0^\i\le( f(x+y) -
f(x) - yf'(x)\mbf 1_{\{|y|<1\}}\ri)\nu(\td y).
$$
Substituting equation (\ref{eq:mkg}) into equation (\ref{ITOPHASE})
and using that $\bff h$ is arbitrary,
we conclude that $\bff Q^+$ solves equation (\ref{systemnoiseWH}).
(ii) Now we turn to the proof of the uniqueness of the
solution of (\ref{systemnoiseWH}). To that end, let ${\bff
G}\in\mc Q(|E|)$ be another solution of (\ref{systemnoiseWH}) and
define the function $\WT f$ by replacing $\bff Q_+$ by ${\bff G}$
in (\ref{eq:mkg}); by an application of It\^{o}'s lemma it follows
then that $\WT f(Y_t, A_t)$ is a local martingale that is bounded
on $[0,\g_k]$ and invoking the optional stopping theorem yields
that
\begin{align}\nn
\WT{f}(j,x) &= \Exp_{j,x}[\WT f(Y_{t\wedge \g_k}, A_{t\wedge\g_k})]\\
&= \Exp_{j,x}[\WT f(\WT Y^+_{k}, A_{\g_k})\mbf 1_{\{\g_k<\i\}}]
+ \lim_{t\to\i}\Exp_{j,x}[\WT f(Y_{t}, A_{t})\mbf 1_{\{\g_k=\i\}}].
\label{optionals}
\end{align}
By definition of $\WT f$ and the absence of positive jumps of $A$,
the first expectation in (\ref{optionals})
is equal to $f(j,x)$.
Note that second term in (\ref{optionals}) is zero if $\bff Q$ is transient or
$\bff Q$ is recurrent and $\sup_t A_t=+\i$. Indeed, in the latter case,
$\g_k$ is finite a.s., whereas in the former case $\Prob(Y_t\in E)$
converges to zero. Thus $f=\WT f$ and the statement (ii) follows.
\exit

\proof{\it of Theorem \ref{uni-s}}
Denote by the 
matrix $\bff Q^+\in\mc Q(|E^+|)$ the generator of $\WT Y^+$
and let $h$ be any function on $E^+$. Then
the Markov property of $\WT Y^+$
yields that $f: E\times\R\to\R$ defined by
$f(i,x) = \Exp_{i,x}[f(\WT Y^+_k)\mbf 1_{\{\g_k<\i\}}]$
is given by
\eqn{
\label{eq:mkgs}
f(i,u) =
\begin{cases}
\bff 1_i'\exp(\bff Q^+(k-x))\bff h & i\in E^+;\\
\bff 1_i'\bff\eta\exp(\bff Q^+(k-x))\bff h & i\in E^-.
\end{cases}
}
Following the line of reasoning of 
the proof of Theorem \ref{uni-ns} and
replacing everywhere (\ref{eq:mkg}) by (\ref{eq:mkgs}),
we arrive at Theorem \ref{uni-s}.\exit

\subsection{Matrix factorisation of $X$}
In case $(A, Y)$ is the embedding as in Section \ref{sec:emb} of a
L\'{e}vy process $X$ of the form (\ref{Levygenal}) we can say more
about the structure of the generator matrix $\bff Q^+$ of the
associated time changed Markov process $\WT Y^+$ in
(\ref{eq:jtilde}). We distinguish between the cases in which
$X^{(-)}$ has non-increasing or non-monotone paths.

In the former case $-X^{(-)}$ is a subordinator and
the additive process $A$ does not increase on the set
$\{s: Y_s=0\}$ and we have $E^- = \{0\}$.
Note that in this case the matrix
$\bff\eta$ in (\ref{etaij})
is of size $1\times(|E|-1)$.
To show explicitly the discounting at rate $a$ we write
\eqn{\label{eq.etaj}
\eta_a(j) = \Exp_{0,0}[\te{-a\g_0}\mbf 1_{\{\WT Y^+_0=j\}}],\q a\ge0,
}
where $\Exp_{0,0}$ denotes the expectation conditioned on $\{A_0=Y_0=0\}$.
As before, for an additive process $A$ of the form
(\ref{eq:AXmin}) to be the embedding of
a L\'{e}vy process $X$ of the form (\ref{Levygenal}),
we set  $\s(i)=\d_{0i}$  and
$v(i)=1-\s(i)$ in (\ref{eq:AXmin}) and let the intensity
matrix $\bff Q_a$ of the special form (\ref{Qa}).
Inserting these quantities in equation (\ref{systemnoiseWHs}) with
$\bff\Sigma=(\d_{0i}, i\in E)_{\mrm{diag}}$
and $\bff V = \bff I - \bff\Sigma$ leads
to the following result:
\begin{Prop} If $-X^{(-)}$ is a subordinator,
the generator matrix $\bff Q_a^+$ of the Markov chain
$\WT Y^+$ in (\ref{eq:jtilde}) corresponding
to the embedding is given by
\eqn{\label{eq:qaplus}
\bff Q_a^+=\bff T^{(+)} + \bff t^{(+)} \bff{\eta}_a,
}
where the vector $\bff\eta_a$ satisfies the equation
\eqn{\label{etaa1}
\bff\eta_a = \l^{(+)}\bff\a^{(+)}\le((\l^{(+)}+a)\bff I -
\k_{X^{(-)}}(-\bff Q_a^+)\ri)^{-1}.
}
\end{Prop}
There is an alternative probabilistic derivation 
of the formula (\ref{eq:qaplus})
for the intensity matrix which uses the Markovian structure.
See \cite{a} for a similar argument in a random walk setting.
Let $\bff M(i,j)$ denote the $ij$th element of a matrix $\bff M$.
By the form of the time change, $\WT Y^+$ is equal to $Y$ 
when $A$ is at a maximum.
Thus, for $i\neq j$,
the total intensity $\bff Q_a^+(i,j)$ of a jump of $\WT Y^+$ from $i$ to $j$,
is given by the intensity $\bff T^+(i,j)$ for the direct 
transitions from $i\to j$
%
added to $(\bff t^{(+)} \bff{\eta}_a')(i,j)$, the intensity
of first passing from $i$ to $\Delta$ and then being renewed
with initial distribution $\bff\eta_a$. The latter is equal to
the  distribution of $\WT Y^+$ at the end of a negative excursion of $A$
away from its supremum. Similarly, the total rate $\bff Q_a^+(i,i)$
of a jump in state $i$ is equal to the rate $\bff T^+(i,i)$ of a direct jump
added to the rate $(\bff t^{(+)} \bff{\eta}_a')(i,i)$ of a jump after being
killed and renewed.

%
In the next section we shall also give a probabilistic argument to prove
the second relation (\ref{etaa1}).  Now we turn to the case
$X^{(-)}$ does {\it not} have decreasing paths. Let $\F(q)$ denote
the largest real root of $\k_{X^{(-)}}(s)=q$, which is positive for $q>0$.

\begin{Prop}\label{prop:Mme}
If $-X^{(-)}$ is not a subordinator,
the generator matrix $\bff Q_a^+$ of the Markov chain
$\WT Y^+$ in (\ref{eq:jtilde}) corresponding
to the embedding is given by
\eqn{\label{qplusam}
\bff Q^+_a = \bff M^{(+)} + \bff m^{(+)} \bff\eta_a
}
where
\eqn{\label{qplusa}
\bff M^{(+)} = \begin{pmatrix}
-\F(a+\l^{(+)}) & \bff 0\\
\bff t^{(+)} & \bff T^{(+)}
\end{pmatrix},
\ \ \bff m^{(+)} = -\bff M^{(+)}\bff 1 = \begin{pmatrix} \F(a+\l^{(+)})\\ 
\bff 0
\end{pmatrix}
}
and $\bff\eta_a$ satisfies
\eqn{\label{etaa21}
\bff{\eta}_a = \frac{\l^{(+)}}{a+\l^{(+)}}
(0,\bff\a^{(+)})\f^-_{X^{(-)}}(-\bff Q^+_a).
}
where
$
\f^-_a(-\bff Q) = \I_0^\i\te{\Bf Q x}\Prob(-I^{(-)}_{e(a)}\in\td x)
$
for generator matrices $\bff Q$
\end{Prop}
\proof The form of the matrices $\bff\Sigma$ and $\bff V$ and
equation (\ref{systemnoiseWH}) imply that all rows of $\bff Q_a^+$
from the second one on till the last one are given (in block
notation) by $(\bff t^{(+)}, \bff T^{(+)})$. Rewriting
(\ref{systemnoiseWH}) we find for the first row
\eqn{\label{whmatrel}
\le[\frac{1}{a+\l^{(+)}}\le(\k_{X^{(-)}}(-\bff Q_a^+) -
(a+\l^{(+)})\bff I\ri)\ri]_1
=-\frac{\l^{(+)}}{a+\l^{(+)}}(0,\bff\a) } where $\bff M_1$ denotes
the first row of a matrix $\bff M$. 
From (\ref{WH2}) with $\bff
T^{(+)}\equiv0$, it follows that $M^{(-)}_{e(q)} = \sup_{s\leq
e(q)}X^{(-)}_s$, the supremum of $X^{(-)}$ at an independent
exponential time $e(q)$, has an exponential distribution with mean
$\F(q)^{-1}$, so that for $q>0$ and $\Re(s)\ge0$
$$
q^{-1} (\k_{X^{(-)}}(s) - q) \f_{X^{(-)}}^-(s) = \F(q)^{-1}(\F(q) - s),
$$
where $\f_{X^{(-)}}^-(s)=\f_{q, X^{(-)}}^-(s)$, 
the moment generating function of the infimum
$I^{(-)}_{e(q)} = \inf_{s\leq e(q)}X^{(-)}_s$ of $X^{(-)}$ at $e(q)$. 
By the Cayley-Hamilton theorem this relation remains valid with
$s$ replaced by a non negative definite matrix
(and thus in particular with $s$ replaced by $-\bff Q_a^+$).
Multiplying both sides of (\ref{whmatrel})
from the right with the matrix
$\f^-_{X^{(-)}}(-\bff Q_a^+)$ yields that
$$
\F(a+\l^{(+)})^{-1}(\F(a+\l^{(+)})\bff I + \bff Q_a^+)_1 =
\frac{\l^{(+)}}{a+\l^{(+)}}(0,\bff\a)\f^-_{a+\l^{(+)}, X^{(-)}}(-\bff Q_a^+).
$$
Thus we find that the first row of $\bff Q_a^+$ is given by
$$
(\bff Q_a^+)_1 = -\F(a+\l^{(+)}) + \F(a+\l^{(+)})
\frac{\l^{(+)}}{a+\l^{(+)}}(0,\bff\a)\f^-_{a+\l^{(+)}, X^{(-)}}(-\bff Q_a^+)
$$
and the result (\ref{qplusam}) -- (\ref{etaa21}) follows.\exit

As above there is also a probabilistic
derivation of the form of $\bff Q_a^+$.
Write $T_1$ for
the last ascending ladder time of $X^{(-)}$ before
the first jump of $X$
and $T_2$ for the first ascending ladder time of $X$ after $T_1$.
and let the vector $\bff\eta_a$ denote the distribution of $\WT Y^+$
at the end of the excursion away from the supremum $X_{T_1} =
M^{(-)}_{T_1}$
$$
\eta_a(k) = \Exp[\te{-a(T_2-T_1)}\mbf 1_{\{Y_{T_2} = k\}}].
$$
%
Note that the supremum of $X$ before the first positive jump of
$X$ (where $X$ is killed at rate $a$ if $a>0$) 
has the same distribution as that of $X^{(-)}$ killed at an 
rate $a + \l^{(+)}$. Since $M^{(-)}$ has an exponential  
with mean $\F(a+\l^{(+)})$, 
the intensity $\mbf Q_a^+(0,i)$, $i\neq 0$, of a jump of
$\WT Y^+$ from $0\to i$ is the intensity $\F(a+\l^{(+)})$
of the jump $0\to\Delta$ times the initial distribution 
$\bff\eta^+_a(i)$ for $\WT Y^+$ 
to be renewed in state $i$ at the end of the excursion of $A$ 
away from its supremum.
As before we see that the rate $\mbf Q_a^+(0,0)$ of $\WT Y^+$ jumping
in state $0$ is the sum of the rate $\F(a+\l^{(+)})$ for $\WT Y^+$
to jump directly and the rate $\F(a+\l^{(+)})\bff\eta_a(0)$
to exit first and be renewed to state $0$.

The other rows follow by the argument
given in the subordinator case. In the next section we shall also give
an alternative derivation of (\ref{etaa21}).

\section{First passage and ladder processes}\label{sec:fluc}
If $X$ is a L\'{e}vy process of the form (\ref{Levygenal}) 
and $X^{(-)}$ is not a subordinator, 
it follows from Proposition \ref{prop:Mme} that the 
running supremum $M_{e(q)}=\sup_{0\leq t\leq e(q)}X_t$ 
of $X$ at an independent exponential time has
a phase-type distribution given by
\eqn{\label{eq:lawmax} \Prob( M_{e(q)} > k) = \bff%
1_0\exp\le\{(\bff M^{(+)} + \bff m^{(+)}\bff\eta_q)k\ri\}\bff 1, }
where $\bff 1_0$ denotes the row-vector with a one in the position corresponding to 
$E^0=\{0\}$ and 0 else. Using alternative proofs, 
this result was found before in Asmussen et al. \cite{AAP}
and Mordecki \cite{MordeckiMax}.
Since $\Prob( M_{e(q)} > k) = \Exp[\te{-q T^+(k)}]$, 
the Laplace transform of the first passage time 
$T^+(k)$, defined in (\ref{eq:fipa}, is equal 
to the left-hand side of (\ref{eq:lawmax}).
 
The result below shows that it is possible to extend
this result to a characterisation and 
description of the law of the up-crossing
ladder process of $X$.
%
%
\begin{Thm}\label{prop:lh}
Let $X$ be a L\'{e}vy process of the form
(\ref{Levygenal}) such that
$-X^{(-)}$ is not a subordinator. Then the following are true:
\BEN
\im[(a)] The ladder height process $H$ is a subordinator given by
\eqn{\label{eq.ladderheight}
H(t) = t +  \sum_{n=1}^{N_t}\WT U_n,
}
where $N$ is a Poisson process with intensity $\F(\l^{(+)})$, 
where $\F(\l^{(+)})$ is the largest root $s$ of $\k_{X^{(-)}}(s) = \l^{(+)}$ and
$\WT U_n$ are i.i.d. random variables with distribution
\eqn{\label{eq.unjumps}
\Prob(\WT U_n\in\td y) = \bff{\eta}_0\cdot (\d_0(\td y),
\exp(\bff T^{(+)}y)\bff t^{(+)}\td y)
}
with $\d_0$ the delta measure in $0$
and $\bff{\eta}_0$ given by (\ref{etaa21}).
\im[(b)] Assume that the representation $(m^{(+)}, \bff T^{(+)}, \bff\a^{(+)})$
is minimal. Then the cumulant $\k^+$ of $(L^{-1},H)$ is given by
\eqn{\label{kappa+}
\k^+(a,s) =  \frac{\prod_{i\in\mc I^{(+)}}(s+\r_i(a))}
{\det(s\bff I-\bff T^{(+)})}
}
for $s$ with $\Re(s)=0$, $a\ge0$,
where the denominator is taken to be one if $\l^{(+)}=0$.
\EEN
\end{Thm}
\begin{Ex}\label{kappa+as}
If there are no positive jumps ($\l^{(+)}=0$),
and $-X^{(-)}$ is not a subordinator, we find back (e.g. \cite{bert96}) that
$$
\k^+(a,s) = s + \F(a)
$$
where $\F(a)$ is the unique positive real root of $\k_{X^{(-)}}(s)=a$.
\end{Ex}

From Theorems \ref{prop:lh} and the extended form of the
Wiener-Hopf factorisation (\ref{WH4}), we can now also determine
the distribution of the downward ladder process by finding an explicit form 
for its cumulant. Inserting the explicit expression (\ref{kappa+}) for $\k^+$
and comparing yields now:

\begin{Cor}
Suppose the representation $(m^{(+)}, \bff\a^{(+)}, \bff T^{(+)})$ is minimal. 
Then the cumulant $\k^-$ of the dual ladder process $(\WH L, \WH H)$
is given by
$$
\k^-(a,s) = C'\times{(a - \k(s))}\frac{\det(s\bff I-\bff T^{(+)})}
{\prod_{i\in\mc I^{(+)}}(s-\r_i(a))}.
$$
where $C'>0$ is some constant.
\end{Cor}
Below we give elementary proof for Theorem \ref{prop:lhs}.

\proof{\it of Theorem \ref{prop:lh}} We first determine the form
of the process $H$. Since $X^{(-)}$ has no positive jumps, $\D M_t
> 0$ implies that $\D X^{(+)}_t>0$ and similarly, if $\td M^c_t > 0$
then $\td M^{(-)}_t>0$ as well, where
$M^{(-)}_t$ denotes the running supremum of $X^{(-)}$ up to time $t$. 
Write $\s_1, \s_2,
\ldots$ for the jump times of $X^{(+)}$ and define recursively 
for $i=1, 2, \ldots$ the stopping times
$G_i = \inf\{s\ge D_{i-1}: M_s=X_s\}$ and $D_i = \inf\{\s_j:
\s_j>G_i\}$, where $G_0=0, D_0 =\s_1$.
Note that, for $t\in[G_i, D_i)$, $M_t$ is continuous and that 
$M$ may jump at $G_i$. 
As $X^{(+)}$ is a compound Poisson process with rate $\l^{(+)}$,
the differences $D_i-G_i$ are exponentially
distributed with parameter $\l^{(+)}$ and, as $X^{(-)}$ has i.i.d. increments 
and $M^{(-)}(e(\l^{(-)}))$ has an exponential distribution with mean $\F(\l^{(+)})^{-1}$,
$M_{D_i} - M_{G_i}$ are i.i.d. exponential with parameter $\F(\l^{(+)})$. 
Since we have taken the local time $L$ to be equal to $M^c$, 
the inter-arrival time of two jumps of $H$ is $\exp(\F(\l^{(+)}))$
distributed and thus the process $H$ is given by
$$
H(t) = M^c(L^{-1}(t)) + \sum_{s\leq t}\D M(L^{-1}(s))\mbf 1_{\{\D M(L^{-1}(s))>0\}}
 = t + \sum_{n=1}^{N_t} \WT U_n
$$
where the $\WT U_n$ are i.i.d. nonnegative random variables (since
$X$ is a L\'{e}vy process)
and $N_t$ an independent Poisson process with rate $\F(\l^{(+)})$.
The jump-size $\WT U_n$ has the same distribution as the
overshoot $X_{T^+(0)}$ of $X$ over the level $0$, 
if $X_0$ is distributed according to $-A+B$ 
where $A$ has the distribution  
$\xi(\td x) = \Prob[(M^{(-)}-X^{(-)})_{e(\l^{(+)})}\in\td x]$ of 
the distance of 
$X^{(-)}$ to its supremum at an exponential time
$e(\l^{(+)})$ and $B$, independent of $A$, 
is distributed according to the jump-sizes of $X^{(+)}$.
Since the upward jumps of $X$ are phase-type, it follows 
that the distribution of the overshoot $X_{T^(0)}$ has an atom in zero and on $(0,\i)$ 
is (defective) phase-type (see e.g. \cite[Prop. 2]{AAP} for a proof).  
A generator matrix of this phase-type distribution is seen to be given by
$\bff T^{(+)}$, with `initial distribution over the phases' 
$\bff\eta_0$ given by the distribution 
of the underlying Markov process at the moment of 
crossing. In the equivalent setting 
of the embedding $(A,Y)=(A_X, Y_X)$ 
of $X$, $\bff\eta_0$ thus satisfies
\begin{align*}
\bff\eta_0(j) &= \Prob[\WT Y^+_0 = j | A_0 \sim \xi, Y_0 \sim (0,\bff\a^{(+)})]\\ 
&= \I_0^\i (0,\bff\a^{(+)})\te{\Bf Q_+ x}\xi(\td x),
\end{align*}
where $j=0, \ldots, m^{(+)}$. Thus the vector $\bff\eta_0$ is given by 
(\ref{etaa21}) and (\ref{eq.unjumps}) is proved.

Finally, we turn to the proof of the identity (\ref{kappa+}). On
the one hand, since $H$ is a compound Poisson process with unit
drift,
$$
\lim_{s\to\i} s^{-1}\k^+(a,s) = 1.
$$
On the other hand, the form of the Wiener-Hopf factor 
$\f_a^+$ and the fact that $|\mc I^{(+)}|=m^{(+)}+1$
imply that $\lim_{s\to-\i}s\f_a^+(s) = (-1)^{m^{(+)}}
\prod_{i\in\mc I^{(+)}}\r_i/\det(\bff T^{(+)})$. 
Combining with (\ref{WH2}) completes the proof.\exit

We end this section with the characterisation 
of the form of the upcrossing ladder if $X^{(-)}$
is the negative of a subordinator.

\begin{Thm}\label{prop:lhs}
Let $X$ be a L\'{e}vy process of the form
(\ref{Levygenal}) such that $X^{(-)}$ is a subordinator.
Then the following are true:
\BEN
\im[(a)] $H$ is a compound Poisson process with jump intensity $c$,
the constant in the definition of $L$,
and jump distribution of phase type with representation
$(\bff\eta^+_0,\bff T^{(+)})$ where $\bff\eta^+_0$
is given by (\ref{etaa21}) with $(0,\bff\a^{(+)})$ replaced by $\bff\a^{(+)}$.
\im[(b)] Suppose the representation $(m^{(+)}, \bff\a^{(+)}, \bff T^{(+)})$ is minimal.
For $s$ with $\Re(s)=0$, $a\ge0$,
\eqn{\label{kappa+s}
\k^+(a,-s) =  c(1 - \bff\eta_a\bff 1)\frac{\prod_{i\in\mc I^{(+)}}(s+\r_i(a))}
{\det(s\bff I-\bff T^{(+)})}\frac{\det(\bff T^{(+)})}{\prod_{i\in\mc I^{(+)}}(\r_i(a))}.
}
where $c>0$ is the constant from the construction of the local time $L$.
\EEN
\end{Thm}

\no\proof{\it of Theorem \ref{prop:lhs}}\ Since in this case
$-X^{(-)}$ is a subordinator, $X$ is irregular for $(0,\i)$, that
is $T = T(0)$, the first entrance time of the positive half line,
is positive a.s. and $X$ will first enter $(0,\i)$ by a jump.
Thus, $(T, X_T)$ are the first ladder time and ladder height
respectively. The independence and homogeneity of the increments
of $X$ imply now that $H$ increases by jumps which sizes are
independent and distributed as $X_T$. By randomisation in the
construction of local time $L$, we deduce that $H$ is a compound
Poisson process with intensity $c$ and jump-sizes distributed as
$X_T$.

We now determine the law of $X_T$.
Denote by $Y$ the underlying Markov process of the jump $\D X_T$.
Then $X_T$ is the lifetime of $Y$
with initial distribution $\bff\eta_0$ given in (\ref{eq.etaj}).
By the defining property of phase type distributions,
$X_T$ is has a distribution that is  phase-type with 
representation $(m^{(+)}, \bff\eta_0,\bff T^{(+)})$.
Since $-X^{(-)}$ is a subordinator independent
of $X^{(+)}$, it holds that $\s$, the first jump time of $X^{(+)}$, 
is exponentially distributed with mean $1/\l^{(+)}$ and
$$
\Prob(X_\s^{(-)}\in\td x, Y_\s=k) =
\I_0^\i\a_k\l^{(-)}\te{-\l^{(-)}t}\Prob(X^{(-)}_t\in\td x)\td t.
$$
Conditioning on the position
of $X^{(-)}$ at time $\s$ and recalling that 
$\WT Y^+$, the upcrossing ladder process of $Y$, is 
a Markov process with generator $\bff Q^+$, it follows that
$$
\bff\eta_0 = \bff\a^{(+)} \I_0^\i\te{\Bf Q_+
x}\Prob(-  X^{(-)}_{e(\l^{(+)})}\in \td x),
$$
where as before $e(q)$ denotes an independent
exponential time with parameter $q$.
Note that the equation for $\bff\eta_0$ is the same as equation (\ref{etaa1})
(with $a=0$) found before.

Finally, we prove the identity (\ref{kappa+s}).
Denote by $G$ the distribution of $T$ and note that in this case
the Laplace transform $\WH G$ of $G$ is $E[\te{-aT}] = \bff\eta_a'\bff 1$.
By the randomisation in the construction of the local time $L$
in this case the Laplace transform of $L^{-1}(1)$ is seen to be
$$
\k^+(a,0) = E[\te{-aL^{-1}(1)}] = \sum_{n=0}^\i P[N_1=n] \WH G^n = \te{c(\WH G-1)},
$$
where $N$ denotes an independent Poisson process with rate $c$
connected to this randomisation. Combining with the factorisation
(\ref{WH2}) and the form (\ref{eq:lawmax}) of the law of the
supremum of $X^{(-)}$ completes the proof in this case.\exit

\section{Nonlinear iteration}\label{sec:iter}
To solve explicitly for the first passage law, the law of the
maximum (\ref{eq:lawmax}) or the law (\ref{eq.unjumps}) of the
up-crossing ladder process of the L\'{e}vy process $X$, we have to
compute the generator matrix $\bff Q_a^+$. One approach achieving this
proceeds by solving the equations (\ref{systemnoiseWHs}) and 
(\ref{systemnoiseWH}) numerically (e.g. via eigenvalue
methods). A different approach,
exploiting the fact that in this case the matrix $\bff Q_a^+$ has a
special structure given by equations (\ref{eq:qaplus}) --
(\ref{qplusam}), is to compute the sub-probability vector
$\bff\eta_a$.

If $X^{(-)}$ is the negative of a subordinator, we consider the
sequence $(\bff\eta^{(n)})_n$ where $\bff\eta^{(0)}$ a
sub-probability vector and $\bff\eta^{(n+1)}$, $n\ge 0$, is
given by the right-hand side of (\ref{etaa1}) with the matrix
$\bff Q^+$ replaced by $\bff T^{(+)} + \bff
t^{(+)}\bff\eta^{(n)}$. If the distribution $F(\td x) = \mbf
1_{\{x\leq 0\}}\Prob(X^{(-)}_{e(a+\l^{(+)}}\in\td x) + \mbf
1_{\{x>0\}}\bff\a\te{\Bf T x}\bff t\td x$ has exponential moments
and non-zero mean, Corollaries 3.2 and 3.3 in  Asmussen \cite{a}
imply $(\bff\eta^{(n)})_n$ converges geometrically fast 
to the solution of (\ref{etaa1}). 
Asmussen \cite{a} proves these results using a coupling argument.

In the sequel we therefore restrict to the case that $X^{(-)}$ is not the
negative of a subordinator. To prove the results we shall follow a route 
that is different from \cite{a}.
Set $\bff M^{(+)}$
equal to the matrix in (\ref{qplusa}) and let $\mc S$ be equal to
the set of sub-probability vectors in $\R^{m+1}$ with $m=m^{(+)}$. 
From 
the interpretation as probability derived in the previous section 
it follows that the function 
$\ps$ given by
\begin{equation}\label{def:psi}
\ps:\bff\eta\mapsto \frac{\l^{(+)}}{a+\l^{(+)}}(0,\bff\a^{(+)})
\f^{(-)}_a(-\bff M^{(+)} - \bff m^{(+)}\bff\eta)
\end{equation}
maps $\mc S$ to $\mc S$.
In the next result it is shown that 
the fixed point equation $\bff\eta = \ps(\bff\eta)$
is uniquely solved by the vector $\bff\eta_a$.
%
%
%
\begin{Thm}\label{thm:iter}Let $a\ge0$.
The following hold true:

(i) The equation $\bff\eta = \ps(\bff\eta)$
has a unique solution $\bff\eta\in\mc S$.

(ii) With $\bff\eta^{(0)}=\mbf 0$ and
$\bff\eta^{(n+1)} = \ps(\bff\eta^{(n)})$,
$\bff\eta^{(n)}\uparrow\bff\eta_a$ as $n\to\i$.
\end{Thm}

\proof (ii) Write $\bff\eta\leq{\bff\eta}'$ if ${\bff\eta}'-\bff\eta$
is non-negative (coordinate-wise).
We claim that $\ps$ satisfies the following monotonicity property: if
$\bff\eta\leq{\bff\eta}'$ then $\ps(\bff\eta) \leq \ps({\bff\eta}')$.
Let $\bff P$ and $\bff P'$ be the 
transition matrices of the Markov chains
with respective generators given by $\bff G(\bff\eta) := \bff M^{(+)} + \bff m^{(+)}\bff\eta$
and $\bff G(\bff\eta')$. Then the matrix
$\bff D = \bff P' - \bff P$ satisfies the matrix differential equation
$$
\dot{\bff D} = \bff D\bff G(\bff\eta) + \bff P'\bff m^{(+)}(\bff\eta'-\bff\eta)\q \bff D(0)=\bff O,
$$
the solution of which is given by 
$\bff D(t) = \I_0^t\bff P'(s)\bff m^{(+)}(\bff\eta'-\bff\eta)\te{(t-s)\Bf G(\Bf \eta)}\td s$.
Hence, coordinate-wise, $\bff D$ is non-negative and 
since 
$$
\psi(\bff\eta) = \frac{\l^{(+)}}{a + \l^{(+)}}(0,\bff\a^{(+)})
\I_0^\i\te{x\Bf G(\bff\eta)}
\Prob(-I^{(-)}_{e(a)}\in\td x)
$$
the claim follows.

Starting with $\bff\eta^{(0)}=\mbf 0$ and setting
$\bff\eta^{(n+1)} = \ps(\bff\eta^{(n)})$ generates a sequence
$(\bff\eta^{(n)})_n$ in $\mc S$ which is nonnegative, since the vector
$\ps(\mbf 0)$ is nonnegative, and coordinate-wise non-decreasing.
Thus, the sequence has a limit in $\mc S$, say $\bff v$, and by
continuity of $\ps$ it follows that $\bff v = \ps(\bff v)$.

(i) Since the matrix $\bff M^{(+)} + \bff m^{(+)}\bff v$, for
$\bff v\in\mc S$, is a generator matrix and solves equation
(\ref{eq:kappamin}), uniqueness follows from Theorem \ref{uni-ns}
if $\k'(0^+)\ge 0$ or $a>0$. Indeed, in the latter case the generator matrix 
$\bff Q_a$ is transient, whereas in the former case $A$ 
satisfies $\sup_t A_t = +\i$.

If $\k'(0^+) < 0$, Theorem \ref{uni-ns} does not apply and we need to provide
a different argument to establish unicity.  
However, using exponential tilting, we shall show that this 
case can be reduced to the case of positive drift.
If $\k'(0^+)<0$, there exists a positive root $\g>0$ of the equation 
$\k(s)=0$. 
Define the tilted measure $\Prob^c$ 
for any positive constant $c\ge \g$ with $\k(c)<\i$ by 
the Radon-Nikodym derivative
$$
\le.\frac{\td\Prob^c}{\td\Prob}\ri|_{\mc F_t} = 
\exp( c X_t - \k(c)t),\q t\ge0,
$$
and denote by $\k^c$, $\f^{-,c}_a$, $\F^c$ and 
$\bff\eta^c$ the respective quantities $\k$, $\f^-_a$, $\F$ and $\bff\eta$
under the measure $\Prob^c$.  It is straightforward to check that, 
for any $c\ge \gamma$, 
$\k^c(s) = \k(s+c)$, $\f^{-,c}(s) = \f^-(s+c)$, $\F^c(a) = \F(a) - c$, 
$\k^{c\prime}(0)=\k'(c)>0$. 
Hence under $\Prob^\g$ the process $X$ has a positive drift and
unicity will follow if we show that  (\ref{etaa21}) can be formulated in terms
of quantities of the process $X$ under the tilted measure $\Prob^\g$.
The next result (from \cite{Asm89})  shows that, under the tilted measure
$\Prob^\g$, the jumps of $X^{(+)}$ remain of phase type:
\begin{Lemma}\label{phasetypetilt} Under $\Prob^\g$ the jumps of $X^{(+)}$ are of phase--type
with representation given by 
$$
(\l^{(+,\g)},\bff\a^{(+,\g)},\bff
T^{(+,\g)}) = (\l^{(+)}\H F^{(+)}[-\g], \bff\a^{(+)}\bff\Delta/\H F^{(+)}[-\g], 
\bff\Delta^{-1}\bff T^{(+)} \bff\Delta + \g\bff I),  
$$
where $\bff\Delta$ is the
diagonal matrix with $k_j$ on the diagonal where 
$\bff k = (\g\bff I-\bff T^{(+)})^{-1}\bff t^{(+)}$. Further, 
$\bff t^{(+,\gamma)}=\bff\D^{-1}\bff t^{(+)}$.
\end{Lemma}
Choosing the killing rate $a^\g = a + \l^{(+)}(1-\hat F^{(+)}[-\gamma])$
and noting that $a + \l = a^\g + \l^\g$ and recalling the form 
of $\bff M^{(+)}$ from (\ref{qplusa}), we deduce that, under $\Prob^\g$, 
the generator matrix of the upcrossing ladder
process $\WT Y^+$ 
is given by $\bff Q^{+,\g} = \bff M^{(+),\g} + \bff m^{(+),\g}\bff\eta_a^\g$, 
where
$$
\bff M^{(+),\g} = \begin{pmatrix} -\F(a+\l^{(+)})+\g&\bff 0\\ 
\bff\D^{-1}\bff t^{(+)} & \bff\D^{-1}\bff T^{(+)}\bff\D + \g\bff I
\end{pmatrix},
$$
and $\bff m^{(+),\g} = ( \F(a+\l^{(+)})+\g, \bff 0)^T$, and $\bff\eta_a^\g$ satisfies 
$\bff\eta = \ps^\g(\bff\eta)$ where $\ps^\g$ is defined as in (\ref{def:psi}) 
but with $a$ replaced by $a^\g$ and all quantities by `tilted' 
ones under the measure $\Prob^\g$. 
Next we provide the link between $\ps$ and $\ps^\g$.
Writing $\bff\D_1 = \text{diag}(1,k_1, \ldots, k_{m^{(+)}})$ and
using Lemma \ref{phasetypetilt}, it follows that
\begin{eqnarray*}
\lefteqn{\ps(\bff\eta)}\\ &=& \frac{\l^{(+)}}{a + \l^{(+)}}
(0,\bff\a^{(+)})
\phi^-_{X^{(-)}}
(-\bff M^{(+)} - \bff m^{(+)}\bff\eta)\\
&=& \frac{\l^{(+)}\H F^{(+)}[-\g]}{a + \l^{(+)}}(0,\bff\a^{(+),\g})\bff\Delta^{-1}_1
\phi^{-,\g}_{X^{(-)}}(-\bff M^{(+)} -
\bff m^{(+)}\bff\eta - \g\bff I)\\
&=& \frac{\l^{(+),\g}}{a + \l^{(+)}}
(0,\bff\a^{(+),\g})\phi^{-,\g}_{X^{(-)}}(-\bff\Delta^{-1}_1
(\bff M^{(+)}+\g\bff I)\bff\Delta_1 -
\bff\Delta^{-1}_1\bff m^{(+)}\bff\eta\bff\Delta_1)\bff\Delta^{-1}_1\\
&=& \frac{\l^{(+),\g}}{a^\g + \l^{(+),\g}}(0,\bff\a^{(+),\g})
\phi^{-,\g}_{X^{(-)}}(
-\bff M^{(+),\g} - \bff m^{(+),\g}\bff{\WT\eta})\bff\Delta^{-1}_1
\end{eqnarray*}
where $\WT{\bff{\eta}} = \bff\eta\bff\Delta_1$.
Since $\k^{\gamma\prime}(0)>0$, Theorem \ref{uni-ns} implies
that $\WT{\bff\eta}=\ps^\g(\WT{\bff\eta})$ has a unique solution 
$\WT{\bff\eta}$ in $\mc S$.
Since it also holds that $\ps^\g(\WT{\bff\eta})$ is equal 
to $\ps(\bff\eta)\Delta_1$ and any solution of 
$\psi(\bff\eta)=\bff\eta$ in $\mc S$ gives rise 
to a solution $\WT{\bff\eta} = \bff\eta\bff\Delta_1$ 
of $\WT{\bff\eta} = \ps^\g(\WT{\bff\eta})$ in $\mc S$,
it follows that, also if $\k'(0)<0$, 
$\bff\eta_a$ is the unique sub-probability vector that solves (\ref{etaa21}). 
\exit

If $a>0$, it turns out that the
convergence of the iteration to its solution
is geometrically fast:

\begin{Thm} Let $a>0$.
If $\k'_{X^{(-)}}(0^+)\ge 0$, the map $\ps$ is a contraction.
If $\k'_{X^{(-)}}(0^+)<0$ and $\k(\F(0))<\i$, 
the map $\ps^{\F(0)}$ is a contraction.
\end{Thm}
\proof
Writing $|\bff\eta| = \sum_i\{|\eta_i|\}$, $\l=\l^{(+)}$, $\bff\a=\bff\a^{(+)}$ 
and $\bff 1_0$ for the column vector that is one in state $0$ and zero else,
it follows by the triangle inequality and the form of $\bff m^{(+)}$ that for
$\bff\eta$, ${\bff\eta}'\in\mc S$
\begin{eqnarray*}
\lefteqn{|\ps(\bff\eta) - \ps({\bff\eta}')|}\\
&\leq& \frac{\l}{a+\l}\I_0^\i
\le|(0,\bff\a)\le[\te{x(\Bf M^{(+)} + \Bf m^{(+)}\Bf\eta)} - 
\te{x(\Bf M^{(+)} + \Bf m^{(+)}{\Bf\eta}')}\ri]\ri|
\Prob(-I^{(-)}_{e(a)}\in\td x)\\
&\leq& \frac{\l}{a+\l}\I_0^\i
\le|(0,\bff\a)\te{x(\Bf M^{(+)} + \Bf m^{(+)}\Bf\xi_x)}
\bff m^{(+)}(\bff\eta - {\bff\eta}')\ri|x\Prob(-I^{(-)}_{e(a)}\in\td x)\\
&\leq& \frac{\l}{a+\l}\I_0^\i
(0,\bff\a)\te{x(\Bf M^{(+)} + \Bf m^{(+)}\Bf\xi_x)}\mbf 1_1
x\Prob(-I^{(-)}_{e(a)}\in\td x)
\F(a)|\bff\eta - {\bff\eta}'|\\
&\leq& \frac{\l}{a+\l}\Exp[-I^{(-)}_{e(a)}]\F(a) |\bff\eta - {\bff\eta}'|,
\end{eqnarray*}
for some vector $\bff\xi_x$ in the convex hull of $\bff\eta$ and
$\bff\eta'$, where we used in the second line the mean value theorem
and in the third line that the integrand is equal to $x$ times the probability
that some Markov chain is at time $x$ in state 1.

The proof is finished noting that $\Exp[-I^{(-)}_{e(a)}]\F(a)\leq 1$
if and only if $\Exp[X^{(-)}_1] \ge 0$. In case $\k_{X^{(-)}}'(0^+)<0$, 
the proof follows by replacing $\ps$ by $\ps^{\F(0)}$ in above reasoning and recalling
that $\k_{X^{(-)}}(\F(0))=0$ and $\k^{\F(0)\prime}_{X^{(-)}}(0)>0$.
\exit

\section*{\bf Acknowledgements.} 
The author would like to thank Florin Avram for inspiring 
and stimulating conversations  and 
Ernst Eberlein for useful comments. Research supported 
by the Nuffield Foundation, grant NAL/00761/G.

\end{document}